\documentclass[12pt,reqno]{amsart}
\newtheorem{theorem}{Theorem}
\newtheorem{proposition}[theorem]{Proposition}

\newtheorem{remark}{Remark}

\newfont{\bb}{msbm10 at 12pt}

\def\pf{{\textit {Proof :} }}

\def\R{\hbox{\bb R}}

\def\Z{\hbox{\bb Z}}
\def\bS{{\mathbb S}}
\def\SM{{\mathbb{S} M}}

\def\Q{{\overline{\mathcal Q}}_\psi}

\newcommand{\bal}{\begin{align}}      \newcommand{\eal}{\end{align}}
\newcommand{\ba}{\begin{array}}      \newcommand{\ea}{\end{array}}
\newcommand{\bc}{\begin{center}}     \newcommand{\ec}{\end{center}}
\newcommand{\be}{\begin{enumerate}}  \newcommand{\ee}{\end{enumerate}}
\newcommand{\beq}{\begin{eqnarray}}  \newcommand{\eeq}{\end{eqnarray}}
\newcommand{\beQ}{\begin{eqnarray*}} \newcommand{\eeQ}{\end{eqnarray*}}
\newcommand{\bi}{\begin{itemize}}    \newcommand{\ei}{\end{itemize}}
\newcommand{\bt}{\begin{tabular}}    \newcommand{\et}{\end{tabular}}
\newcommand{\bdm}{\begin{displaymath}} \newcommand{\edm}{\end{displaymath}}

\def\qed{\hfill{Q.E.D.}\smallskip}

\begin{document}

\title{Dirac Operator on Embedded Hypersurfaces}

\author{Oussama Hijazi}
\address[Hijazi]{Institut {\'E}lie Cartan\\
Universit{\'e} Henri Poincar{\'e}, Nancy I\\
B.P. 239\\
 54506 Vand\oe uvre-L{\`e}s-Nancy Cedex, France}
\email{hijazi@iecn.u-nancy.fr}

\author{Sebasti{\'a}n Montiel}
\address[Montiel]{Departamento de Geometr{\'\i}a y Topolog{\'\i}a\\
Universidad de Granada\\
18071 Granada \\
Spain}
\email{smontiel@goliat.ugr.es}

\author{Xiao Zhang}
\address[Zhang]{Institute of Mathematics\\
Academy of
Mathematics and Systems Sciences, Chinese Academy of
Sciences\\
Beijing 100080, P.R. China}
\email{xzhang@math08.math.ac.cn}

\begin{abstract}
New extrinsic lower bounds are given for the classical Dirac operator
on the boundary of a compact domain of a spin manifold. The main tool is to 
solve some boundary problems for the Dirac operator of the domain under boundary
conditions of Atiyah-Patodi-Singer type. Spinorial techniques are used to give 
simple proofs of classical results for compact embedded hypersurfaces.
\end{abstract}



\thanks{Research of S.M. is partially 
supported by a DGICYT grant No. PB97-0785. Research of X.Z. is
partially supported by the Chinese NSF and mathematical physics
program of CAS}

\thanks{This work is partially done during the visit of the last two
authors to the Institut {\'E}lie Cartan, Universit{\'e} Henri
Poincar{\'e}, Nancy 1.  They would like to thank the institute
for its hospitality.}


\maketitle \pagenumbering{arabic}
 
\section{Introduction}

The spectrum of the fundamental Dirac operator on closed manifolds
have been extensively studied over the past three decades. First, the 
intrinsic aspect has been systematically studied by many authors 
(see \cite{BFGK,BHMM} and for references therein). The striking fact 
in this setup is the Lichnerowicz' Theorem which gives topological
obstructions to the existence of metrics with positive scalar curvature.
Another important feature in this approach is the geometric 
characterisation of manifolds
admitting solutions of some special field equations as the Killing spinor 
equation. 

Second, the extrinsic aspect has been recently studied in
\cite{An,Bm,Ba} where mainly extrinsic upper bounds are obtained.

More recently in \cite{Z,HZ,Mo}, extrinsic lower bounds 
for the hypersurface Dirac operator are established.

In this paper, we investigate the spectral properties of the Dirac 
operator on a compact  manifold with boundary for the 
Atiyah-Patodi-Singer type boundary condition corresponding to the spectral 
resolution of the classical Dirac operator of the boundary hypersurface.  
We start by recalling
the  Schr{\"o}dinger-Lichnerowicz' integral formula (\ref{boun-weit-twis-iden}) 
for a compact $(n+1)$--dimensional manifold 
$\Omega$ with boundary $\Sigma=\partial\Omega$ from which we deduce 
a spinorial Reilly type inequalities (\ref{boun-weit-twis-ineq}) and
(\ref{boun-weit-emt-ineq}).  Under 
some curvature assumptions, we show that the  
Dirac operator on the ambient space $\Omega$, subject to an APS type
boundary condition, has zero kernel and
we derive the analogue of the Friedrich inequality (\ref{frie-ineq})
with its generalization (\ref{hija-ineq}) involving the energy-momentum 
tensor.

We then use the Reilly type Inequality (\ref{boun-weit-twis-ineq}) to
prove that (see Theorem \ref{infH}) if the scalar curvature of $\Omega$ 
and the mean curvature
$H$ of $\Sigma$ are nonnegative, then the lowest nonnegative eigenvalue
of the intrinsic hypersurface Dirac operator is at least equal to 
$$
\frac{n}{2}\inf_\Sigma H.
$$  
It is shown that 
this estimate improves previous results. In particular, it is valid
even if the scalar curvature of the boundary is negative.

Finally, we make use of the spinorial techniques to study constant mean 
curvature and minimal embedded hypersurfaces. A spinorial simple proof 
for the classical Alexandrov Theorem (see Theorem \ref{theo-alex}) is given.
We then prove that (see Theorem \ref{theo-tota-geod}) minimal
compact hypersurfaces bounding a compact domain admitting a parallel
spinor are totally geodesic. We end with a rigidity type result.

The conformal aspect of the Dirac operator on embedded hypersurfaces
is the object of \cite{HMZ1} while  general lower bound estimates 
on compact manifolds with boundary are established in \cite{HMZ2}.

\section{Preliminaries} 
Let $(M, \langle\;,\;\rangle)$ be an $(n+1)$-dimensional Riemannian 
spin manifold and  denote by  $\overline{\nabla}$ its  Levi-Civit{\`a} 
connection. We use the same symbol to denote the corresponding lift 
of $\overline{\nabla}$ to the spinor bundle $\SM$ of $M$. On the spinor 
bundle $\SM$ there exist a natural Hermitian structure
(also denoted by $\langle\;,\;\rangle$), and a Clifford module 
structure\break $\gamma:\hbox{Clif}(M)\rightarrow \hbox{End}(\SM)$ 
which are compatible with $\overline{\nabla}$. That is
\begin{eqnarray}
X\langle\psi,\varphi\rangle &=& \langle\overline{\nabla}_X\psi,
\varphi\rangle+\langle\psi,\overline{\nabla}_X\varphi
\rangle\label{comp1}\\  
\langle\gamma(X)\psi,\gamma(X)\varphi\rangle &=& |X|^2
\langle\psi,\varphi\rangle\label{comp2}\\
\overline{\nabla}_X\Big(\gamma(Y)\psi\Big) 
&=& \gamma(\overline{\nabla}_XY)\psi
+\gamma(Y)\overline{\nabla}_X\psi,\label{comp3}
\end{eqnarray}
for any spinor fields $\psi,\varphi\in\Gamma(\SM)$ and any 
tangent vector fields\break
$X,Y\in \Gamma (TM)$.
The Dirac operator $\overline{D}$ on $\SM$ is locally given by 
\beq
\overline{D}=\sum_{i=1}^{n+1}\gamma(e_i)\overline{\nabla}_{e_i},
\eeq
 where $\{e_1,\dots,e_{n+1}\}$ is a local orthonormal 
frame of $TM$.
Consider  an orientable hypersurface $\Sigma$ in $M$. 
Let $\nabla$ be the Levi-Civit{\`a} connection of the Riemannian 
metric on $\Sigma$ induced by the metric of $M$. The Gauss 
formula says that
\beq\label{riem-gaus}
\nabla_XY=\overline{\nabla}_XY-\langle AX,Y\rangle N,
\eeq
where $X,Y$ are vector fields tangent to the hypersurface 
$\Sigma$, the vector field $N$ is the global unit field (inner) 
normal to 
$\Sigma$ and $A$ stands for the shape operator corresponding 
to $N$, that is,
\beq\label{shap-oper}
\overline{\nabla}_XN=-A X,\qquad \forall X\in \Gamma(T\Sigma).
\eeq
The spin structure of $M$ can be also induced on $\Sigma$ in 
such a way that the restricted bundle $\SM_{|\Sigma}$ is 
isomorphic to either $\bS\Sigma$ or $\bS\Sigma\oplus \bS\Sigma$ 
according to  the dimension $n$ of $\Sigma$ is either even or odd
(see \cite{Ba}, \cite{Mo} for example). 
{\it A spinor field on $M$ and its restriction to the hypersurface 
will be denoted by the same symbol}.
Since the $n$--dimensional Clifford algebra is  the even part of 
the $(n+1)$--dimensional Clifford algebra, Clifford multiplication 
on $\SM_{|\Sigma}$ is given by 
$$
\gamma^\Sigma(X)\psi=\gamma(X)\gamma(N)\psi,
$$
where $\psi\in \Gamma(\SM_{|\Sigma})$ and $X\in \Gamma(T\Sigma)$. 
It is not difficult to check  
that for any $X\in\Gamma(T\Sigma)$ and $\psi\in\Gamma(\SM_{|\Sigma})$, 
the Levi-Civit{\`a} connection on 
$\bS \Sigma$ is given by the following spinorial Gauss formula
\begin{equation}\label{spin-gaus}
\nabla_X\psi=\overline{\nabla}_X\psi-\frac{1}{2}\gamma^\Sigma(A X)\psi
=\overline{\nabla}_X\psi-
\frac{1}{2}\gamma(A X)\gamma(N)\psi\,.
\end{equation}
Hence, if $D$ denotes the Dirac operator associated with the spin 
structure of the hypersurface $\Sigma$, then for any 
spinor field $\psi\in\Gamma(\SM_{|\Sigma})$
\begin{equation}\label{dira-extr1}
D\psi=\sum_{j=1}^n\gamma^\Sigma(e_j)\nabla_{e_j}\psi
=\frac{n}{2}H\psi-\gamma(N)\sum_{j=1}^n\gamma(e_j)
\overline{\nabla}_{e_j}\psi,
\end{equation}
where $\{e_1,\dots,e_n\}$ is a local orthonormal frame on $\Sigma$ 
such that \break$\{e_1,\dots,e_n, e_{n+1}=N\}$
is the corresponding local orthonormal frame on $M$ and 
$$
H=\frac{1}{n}\hbox{trace\,}A$$
is the mean curvature of $\Sigma$ corresponding to the 
orientation $N$. From (\ref{dira-extr1}), if 
$\psi\in \Gamma(\SM)$ is a spinor field on the ambient 
manifold $M$, we have
\begin{equation}\label{dira-extr2}
D\psi=\frac{n}{2}H\psi-\gamma(N)\overline{D}\psi
-\overline{\nabla}_N\psi.
\end{equation}

We end this section by showing that the sepctrum of the Dirac operator
$D$ on the hypersurface $\Sigma$ is symmetric w.r.t zero. We have

\begin{proposition} For any spinor field $\psi\in \Gamma(\SM)$, and 
any tangent vector field $X\in\Gamma(T\Sigma)$, the 
following relations hold
{\setlength\arraycolsep{2pt}
\beq\label{comp-nabl}
\nabla_X\,\big(\gamma(N) \psi\big) 
& = & \gamma(N)\,\nabla_X \psi,\\
    D\,\big(\gamma(N) \psi\big)\label{anti-comm}
& = & -\gamma(N)\,D \psi.
\eeq}
\end{proposition}

\pf By (\ref{spin-gaus}) and (\ref{shap-oper}) it follows
\beQ
\nabla_X\,\big(\gamma(N) \psi\big) 
& = & \Big(\overline{\nabla}_X -
\frac{1}{2}\gamma(A X)\gamma(N)\Big) \big(\gamma(N)\psi\big)\\
& = & -\gamma(A X) \psi + \gamma(N) \overline{\nabla}_X \psi 
+\frac{1}{2}\gamma(A X)\psi \\
& = & \gamma(N) \Big( \overline{\nabla}_X  -
\frac{1}{2}\gamma(A X)\gamma(N)\Big)\psi\\
& = & \gamma(N) \nabla_X \psi.
\eeQ
For the second relation it is sufficient to use (\ref{comp-nabl}) and the 
Clifford algebra relations to get
\beQ
D\,\big(\gamma(N) \psi\big) & = & \sum_{j=1}^n\gamma(e_j)\gamma(N)\nabla_{e_j}
\big( \gamma(N) \psi\big)\\
& = & \sum_{j=1}^n\gamma(e_j)\gamma(N)\gamma(N) \nabla_{e_j}\psi\\
& = & - \gamma(N)\sum_{j=1}^n\gamma(e_j)\gamma(N) \nabla_{e_j}\psi\\
& = & - \gamma(N) D\psi.
\eeQ 

\qed 

\section{Bounding domains hypersurfaces}
Suppose now that the hypersurface $\Sigma$ is the boundary of a 
compact domain $\Omega$ in the manifold $M$ (which could be the 
manifold itself). Let $\psi\in\Gamma(\bS\Omega)$ be a spinor 
field on the domain $\Omega$. The Schr{\"o}dinger--Lichnerowicz 
formula says that
\beq\label{schr-lich}
\overline{D}^2\psi=\overline{\nabla}^*\overline{\nabla}\psi
+\frac{1}{4}\overline{R}\psi,
\eeq
where $\overline R$ is the scalar curvature of $M$. Then$$
\langle\overline{D}^2\psi,\psi\rangle=\langle\overline{\nabla}^* 
\overline{\nabla}\psi,\psi
\rangle+\frac{1}{4}\overline{R}|\psi|^2.$$
Consider the 1-forms $\alpha$ and $\beta$ on $\Omega$ defined by
$$
\alpha(X)=\langle\gamma(X)\overline{D}\psi,\psi\rangle,
\qquad \beta(X)=\langle\overline{\nabla}_X\psi,\psi\rangle,
$$
for any $X\in \Gamma(T\Omega)$. It is clear that
$$
\delta\alpha=\langle\overline{D}^2\psi,\psi\rangle
-|\overline{D}\psi|^2,
$$
because $\gamma$ acts by skew--symmetric endomorphisms for 
$\langle\;,\;\rangle$. Also$$
\delta\beta=-\langle\overline{\nabla}^* \overline{\nabla}\psi
,\psi\rangle+|\overline{\nabla}
\psi|^2.$$
Hence we obtain
$$
\delta\alpha+\delta\beta=|\overline{\nabla}\psi|^2
-|\overline{D}\psi|^2+\frac{1}{4}\overline{R}|\psi|^2.
$$
Integrating on $\Omega$ and applying the divergence theorem
\beQ
-\int_{\Sigma}\langle\gamma(N)\overline{D}\psi + \overline{\nabla}_N\psi,\psi\rangle\,d\Sigma 
=\int_\Omega\Big(
|\overline{\nabla}\psi|^2-|\overline{D}\psi|^2+\frac{1}{4}\overline{R}
|\psi|^2\Big)\,d\Omega,
\eeQ
where $N$ is the inner unit normal field along $\Sigma$.
Using (\ref{dira-extr2}), this equation could be written as
\beq
\label{boun-weit}
\int_\Sigma\Big(\langle D\psi,\psi\rangle
-\frac{nH}{2}|\psi|^2\Big)\,d\Sigma 
=\int_\Omega\Big(|\overline{\nabla}\psi|^2
-|\overline{D}\psi|^2
+\frac{1}{4}\overline{R}|\psi|^2\Big)\,d\Omega,
\eeq
for any spinor field $\psi\in\Gamma(\bS\Omega)$.

On the other hand, for any spinor field $\psi$ on $M$, we 
have the following 
decomposition :
\begin{equation}\label{deco} |\overline{\nabla}\psi|^2 
= |\overline{P}\psi|^2
+\frac{1}{n+1}|\overline{D}\psi|^2, 
\end{equation}
where $\overline{P}$ is the Twistor operator of $M$ 
defined by
\begin{equation}\label{twis-oper}
\overline{P}_X\psi := \overline{\nabla}_X\psi
+\frac{1}{n+1}\gamma(X)\overline{D}\psi,\qquad \forall X\in \Gamma(TM).
\end{equation}
A non-trivial spinor field $\psi$ such that 
$\overline{P}\psi\equiv 0$ is {\it called} a twistor-spinor. 
Combining the identities (\ref{boun-weit}) and 
(\ref{deco}), it follows
\beq\label{boun-weit-twis-iden}
\int_\Sigma\left(\langle D\psi,\psi\rangle
-\frac{nH}{2}|\psi|^2\right)\,d\Sigma
 &= &\frac{1}{4}\int_\Omega \overline{R}|\psi|^2\,
d\Omega \\ & & 
-\frac{n}{n+1}\int_\Omega|\overline{D}\psi|^2\,d\Omega 
+ \int_\Omega|\overline{P}\psi|^2\,d\Omega,
\nonumber
\eeq
for all $\psi\in\Gamma(\bS\Omega)$.

\begin{remark} Since $|\overline{P}\psi|^2 \ge 0$, identity 
(\ref{boun-weit-twis-iden}) immediately translates to
\begin{equation}\label{boun-weit-twis-ineq}
\int_\Sigma\left(\langle D\psi,\psi\rangle
-\frac{nH}{2}|\psi|^2\right)\,d\Sigma\ge \frac{1}{4}
\int_\Omega \overline{R}|\psi|^2\,d\Omega-\frac{n}{n+1}
\int_\Omega|\overline{D}\psi|^2\,d\Omega,
\end{equation}
which is the analogue of the Reilly Inequality \cite{Re} for the 
gradient of a function. Moreover, equality occurs if and only 
if $\psi$ is a 
twistor--spinor.
\end{remark} 

We now make use of the energy-momentum tensor to derive another useful
expression of the r.h.s of identity (\ref{boun-weit}) (see \cite{Hi}). Recall
that the energy-momentum tensor $\Q$ associated with a spinor field
$\psi \in\Gamma(\bS\Omega)$ is the symmetric $2$-tensor,  defined on the 
complement set of zeros of $\psi$ and for any tangent vector fields 
$X, Y\in \Gamma (T \Omega)$ by
\beq\label{ener-mome-tens}
\Q (X , Y) = \frac12 \Re \;\langle\, \gamma(X) \overline{\nabla}_Y \psi 
+ \gamma(Y) \overline{\nabla}_X \psi , \frac{\psi}{\vert\psi\vert^2} \,\rangle.
\eeq
If the associated symmetric endomorphism of the tangent bundle $T \Omega$
is denoted by the same symbol, then one can easily check that the modified 
connection 
$\overline{\nabla}^{\Q}$ defined by
\beq
\overline{\nabla}_X^{\Q}\, \psi : = \overline{\nabla}_X\psi + \gamma(\Q(X)) \psi
\eeq
satisfies, for any spinor field $\psi$ the relation
\beq\label{ener-mome-tens-deco}
\vert\overline{\nabla} \psi\vert^2 = \vert\overline{\nabla}^{\Q}\, \psi\vert^2
+ \vert\Q \vert^2 \vert\psi \vert^2.
\eeq
Hence, this identity combined with (\ref{boun-weit}) imply
\beq
\label{boun-weit-emt}
\!\!\!\!\!\!\!\!\!\!\!\!\!\!\!\!\!\!\!\int_\Sigma\Big(\langle D\psi,\psi\rangle
-\frac{nH}{2}|\psi|^2\Big)\,d\Sigma 
& = & \nonumber\\
& &\!\!\!\!\!\!\!\!\!\!\!\!\!\!\!\!\!\!\!\!\!\!\!\!\!\!\!\!\!\!\!\!\!\!\!\!
\!\!\!\!\!\!\!\!\!\!\!\!\!\!\!\!\!\!\!\!\!\!\!\!\!\!
\int_\Omega\Big(
\big(\frac{1}{4}\overline{R}+ \vert\Q \vert^2\big) |\psi|^2
-|\overline{D}\psi|^2
+|\overline{\nabla}^{\Q}\,\psi|^2\Big)\,d\Omega.
\eeq

\begin{remark} As before, since $|\overline{\nabla}^{\Q}\,\psi|^2\ge 0$, identity 
(\ref{boun-weit-emt}) translates to
\beq\label{boun-weit-emt-ineq}
\!\!\!\!\!\!\!\!\!\!\!\!\!\!\!\!\!\!\!\int_\Sigma\Big(\langle D\psi,\psi\rangle
-\frac{nH}{2}|\psi|^2\Big)\,d\Sigma 
& \ge & \nonumber\\
& &\!\!\!\!\!\!\!\!\!\!\!\!\!\!\!\!\!\!\!\!\!\!\!\!\!\!\!\!\!\!\!\!\!\!\!\!
\!\!\!\!\!\!\!\!\!\!\!\!\!\!\!\!\!\!\!\!\!\!\!\!\!\!
\int_\Omega\Big(
\big(\frac{1}{4}\overline{R}+ \vert\Q \vert^2\big) |\psi|^2
-|\overline{D}\psi|^2\Big)\,d\Omega.
\eeq
Moreover, if equality occurs in  (\ref{boun-weit-emt-ineq}) for a spinor
field $\psi$, then the function $|\psi|^2$ is constant.
\end{remark}

\section{Boundary Problems for the Dirac operator}

Since the hypersurface $\Sigma=\partial\Omega$ is compact, 
the Dirac operator $D$ has a discrete spectrum
$$
\cdots\le\lambda_{-k}\le\cdots\le\lambda_{-1}\le 
0\le\lambda_1\le\cdots\le\lambda_k\le\cdots.
$$
But from  (\ref{anti-comm}), 
we have  $\lambda_{-k}=-\lambda_k$ for all $k\in
\Z$. That is, the spectrum of $D$ is 
$$
\cdots\le-\lambda_k\le\cdots\le-\lambda_1\le 
0\le\lambda_1\le\cdots\le\lambda_k\le\cdots \nearrow.
$$
Denote by $\pi_+:\Gamma(\bS\Sigma)\rightarrow\Gamma(\bS\Sigma)$ 
the projection onto the subspace of 
$\Gamma(\bS\Sigma)$ spanned by the eigenspinors corresponding 
to the nonnegative eigenvalues of $D$. It is clear that
\begin{equation}\label{proj}
D\pi_+=\pi_+D\qquad\hbox{and}\qquad \int_\Sigma\langle D\psi,
\psi\rangle\,d\Sigma\le \int_\Sigma\langle D\pi_+\psi,\pi_
+\psi\rangle\,d\Sigma,
\end{equation}
for any spinor field $\psi$ on $\Sigma$ and the equality holds 
if and only if $\pi_+\psi=\psi$. This projection $\pi_+$ provides 
an Atiyah-Patodi-Singer type boundary condition for the Dirac operator
$\overline{D}$ of the domain $\Omega$.   We have proved 
in \cite{HMZ2} that 
this is a global self-adjoint elliptic condition. This can be done either 
by using standard facts on pseudo-differential operators \cite{BW,GLP,S} 
or by obtaining basic elliptic estimates and standard results from functional
analysis. The second approach was discovered 
in \cite{FS} and closely followed in \cite{HMZ2}. Following either one
of these approaches, we prove the following result:

\begin{theorem}\label{theo-FS} 
Let $\Omega$ be a compact Riemannian spin manifold with boun\-dary 
$\partial\Omega=\Sigma.$ The inhomogeneous boundary problem for 
the Dirac operator
\beq\label{inho-prob}
\left\{
\begin{array}{lll}
\overline{D}\psi&=\Psi \qquad &\hbox{ {\rm on} }\Omega\\
\pi_+\psi&=\pi_+\varphi\qquad &\hbox{ {\rm on} }\Sigma,
\end{array}
\right. 
\eeq
has a smooth solution for each $\Psi\in\Gamma(\bS\Omega)$ and 
$\varphi\in\Gamma(\bS\Omega_{|\Sigma})$ satisfying the following 
integrability condition
$$
\int_\Omega\langle\Psi,\Phi\rangle\,d\Omega
+\int_\Sigma\langle\gamma(N)\varphi,\Phi\rangle\,
d\Sigma=0
$$
for any harmonic spinor $\Phi$ on $\Omega$ such that $\pi_+\Phi=0$. 
This solution is unique up to an arbitrary harmonic spinor field 
$\Phi$ of this type.
\end{theorem}

Under some curvature assumptions, we now show the following particular 
case of Theorem \ref{theo-FS}:

\begin{proposition}\label{prop-APS-boun}
Let $\Omega$ be a compact Riemannian spin manifold with nonnegative 
scalar curvature $\overline{R}$, whose boundary 
$\partial\Omega=\Sigma$ has nonnegative mean curvature $H$
(with respect to the inner normal). Then the following 
inhomogeneous problem  for the Dirac operator $\overline{D}$ 
of $\Omega$ with the Atiyah-Patodi-Singer boundary condition
\beq
\left\{
\begin{array}{lll}
\overline{D}\psi&=\Psi \qquad &\hbox{ {\rm on} }\Omega\\
\pi_+\psi&=\pi_+\varphi\qquad &\hbox{ {\rm on} }\Sigma,
\end{array}
\right. 
\eeq
has a unique smooth solution for any $\Psi\in\Gamma(\bS\Omega)$ 
and $\varphi\in\Gamma(\bS\Omega_{|\Sigma})$.
\end{proposition}

\pf
Since $\overline{R}$ and $H$ are both nonnegative, Inequality 
(\ref{boun-weit-twis-ineq})  could be written as
$$
\int_\Sigma\langle D\psi,\psi\rangle\,d\Sigma\ge 
-\frac{n}{n+1}\int_\Omega|\overline{D}
\psi|^2\,d\Omega,$$
for any spinor field $\psi\in\Gamma(\bS\Omega)$. Assume that $\overline{D}\psi=0$ and $\pi_+\psi=0$. Then, 
$$
0=\int_\Sigma\langle D\pi_+\psi,\pi_
+\psi\rangle\,d\Sigma \ge \int_\Sigma\langle D\psi,\psi\rangle\,d\Sigma
\ge-\frac{n}{n+1}\int_\Omega|\overline{D}
\psi|^2\,d\Omega =0.
$$
Hence, since we have equality in (\ref{boun-weit-twis-ineq}), 
$\psi$ is a harmonic twistor--spinor, that is a parallel spinor. In
particular the function $|\psi|^2$ is constant.
On the other hand, since the equality on the right side of (\ref{proj}) is 
achieved, we have  $\psi=\pi_+\psi=0$ on $\Sigma$.  
Therefore $\psi$ is identically zero on $\Omega$ and the integrability 
condition of Theorem \ref{theo-FS} is satisfied.

\qed

\section{Extrinsic Lower Bounds for the 
Hypersurface Dirac Operator}

We start this section by proving a Friedrich type inequality \cite{Fr} 
for compact spin manifolds with non-empty boundary.

\begin{theorem}\label{theo-frie-boun}
Let $\Omega$ be a compact $(n+1)$--dimensional Riemannian spin  manifold 
of nonnegative scalar curvature $\overline{R}$
with  boundary $\partial\Omega$ of nonnegative mean curvature. Then the 
first eigenvalue $\overline{\lambda}_1$ 
of the Dirac operator on $\Omega$, with the
Atiyah-Patodi-Singer condition, satisfies
\beq\label{frie-ineq}
\overline{\lambda}_1^2 > \frac{n+1}{4n}\inf_\Omega \overline{R}.
\eeq
\end{theorem}

\pf
Consider the Dirac operator $\overline{D}$ acting on the 
spinor bundle $\Gamma(\bS\Omega)$ satisfying the Atiyah-Patodi-Singer 
boundary condition
$$
\pi_+\psi=0\qquad\hbox{on }\partial\Omega=\Sigma.
$$
The spectrum of $\overline{D}$ consists of entirely 
isolated real eigenvalues with finite multiplicity and 
smooth eigenspinors (see \cite{BW,FS} or 1.5.8 in \cite{GLP}). 
Let $\overline{\lambda}_1$ be the eigenvalue of $\overline{D}$ 
with the lowest absolute value and take a corresponding eigenspinor $\psi$. 
Then  from (\ref{proj}), we have 
$$
\overline{D}\psi=\overline{\lambda}_1\psi,\qquad \int_\Sigma\langle 
D\psi,\psi\rangle\,
d\Sigma\le 0.
$$
Inequality (\ref{boun-weit-twis-ineq}) applied to this particular spinor 
field $\psi$ could be written as
$$
0\ge -\frac{n}{2}\int_\Sigma H|\psi|^2\,d\Sigma\ge \frac{1}{4}
\int_\Omega \overline{R} |\psi|^2\,d\Omega
-\frac{n}{n+1}\overline{\lambda}_1^2\int_
\Omega |\psi|^2\,d\Omega.
$$
Hence
\beq\label{frie-boun}
\overline{\lambda}_1^2 \ge \frac{n+1}{4n}\inf_\Omega \overline{R}.
\eeq
If equality in (\ref{frie-boun}) is achieved, then $\psi$ is a Killing 
spinor since it is simultaneously a twistor--spinor 
and an eigenspinor. Moreover, $\psi=\pi_+\psi$ and $H=0$. Since a 
 (real) Killing spinor has constant length and $\psi=\pi_+\psi=0$ on
the boundary $\Sigma$, the spinor field $\psi\equiv 0$ on $\Omega$. This 
contradicts the fact that $\psi$ is a non-trivial eigenspinor of 
$\overline D$. Therefore,
equality in (\ref{frie-boun}) could not hold.

\qed 

The next result is an immediate generalization of (\ref{frie-ineq}) whose
proof follows the same arguments as Theorem \ref{theo-frie-boun}, when
inequality (\ref{boun-weit-emt-ineq}) is used instead of 
(\ref{boun-weit-twis-ineq}). More precisely, we have

\begin{theorem}\label{theo-emt-boun}
Let $\Omega$ be a compact $(n+1)$--dimensional Riemannian spin manifold 
with scalar curvature $\overline{R}$ satisfying 
$\frac{1}{4}\overline{R}+ \vert\Q \vert^2 \ge 0$
whose boundary $\partial\Omega$ is of nonnegative mean curvature. Then the 
first eigenvalue $\overline{\lambda}_1$ 
of the Dirac operator on $\Omega$, with the
Atiyah-Patodi-Singer condition, satisfies
\beq\label{hija-ineq}
\overline{\lambda}_1^2 > \inf_{\Omega_\psi} \,\Big( \frac{1}{4}\overline{R}+ \vert\Q \vert^2 \Big),
\eeq
where $\Omega_\psi$ is the complement of the set of zeros of the associated
eigenspinor field $\psi$.
\end{theorem}

Now, we use the Reilly type inequality 
(\ref{boun-weit-twis-ineq}) to get a lower bound for the first
eigenvalue of the Dirac operator 
$D$ on the boundary hypersurface $\Sigma$ of the compact 
manifold $\Omega$. More precisely, we have

\begin{theorem}\label{infH} Let $M$ be a Riemannian spin manifold
of nonnegative scalar curvature $\overline{R}$ and 
$\Sigma$ a compact hypersurface.
Suppose that $\Sigma$ has nonnegative mean curvature $H$ with 
respect to its inner unit normal field $N$ and that it bounds a 
compact domain $\Omega$ in $M$. Then, the lowest nonnegative 
eigenvalue $\lambda_1$ of the Dirac operator on $\Sigma$ 
satisfies
\beq\label{lowe-boun-H}
\lambda_1\ge\frac{n}{2}\inf_\Sigma H.
\eeq
Moreover, 
if the equality holds, then $\Omega$ is a Ricci flat manifold, 
$\Sigma$ has constant mean curvature and the 
eigenspace corresponding to $\lambda_1$ consists of the 
restrictions to $\Sigma$ of  parallel spinors on the 
domain $\Omega$.
\end{theorem}
\pf {}From Theorem \ref{theo-frie-boun} above or directly 
from Proposition \ref{prop-APS-boun}, we have
$\overline{\lambda}_1 >0$ 
and so the following inhomogeneous boundary problem has a 
unique solution:
\beq\label{inho-boun-cond}
\left\{
\begin{array}{rll}
\overline{D}\psi&=0\qquad &\hbox{on }\Omega\\
\pi_+\psi&=\pi_+\varphi=\varphi\qquad &\hbox{on }
\partial\Omega=\Sigma,
\end{array}
\right. 
\eeq
where $\varphi\in\Gamma(\bS\Sigma)$ is an eigenspinor on 
$\Sigma$ corresponding to the first eigenvalue 
$\lambda_1\ge 0$ of $D$. That is
$$
D\varphi=\lambda_1\varphi\qquad\hbox{and so}
\quad \pi_+\varphi=\varphi.
$$
{}From the Reilly inequality (\ref{boun-weit-twis-ineq}), 
we get
$$
\int_\Sigma\left(\lambda_1-\frac{n}{2}H\right)
|\psi|^2\,d\Sigma\ge
\frac{1}{4}\int_\Omega \overline{R}|\psi|^2\,d\Omega,
$$
which implies (\ref{lowe-boun-H}). For the equality case in 
(\ref{lowe-boun-H}) is satisfied for
 a spinor field $\psi$, then  $\psi$ is
harmonic spinor and a twistor--spinor, hence parallel. As 
$\pi_+\psi=\varphi$ along 
the boundary $\Sigma$, $\psi$ is a non-trivial parallel spinor and 
also $\lambda_1=nH/2$. 
The existence of such a spinor field implies that $\Omega$ is a 
Ricci flat Riemannian manifold (see, for example, \cite{Ht} 
or \cite{W}). On the other hand, since $\psi$ is parallel, 
one deduces from (\ref{dira-extr2}) that 
$D\psi=(nH/2)\psi$. Hence, as the equality $\lambda_1=nH/2$
implies that $H$ is constant, we have 
$$
\varphi=\pi_+\psi=\psi.
$$
Conversely, the fact that the restriction to $\Sigma$ of 
a parallel spinor on $\Omega$ is an eigenspinor with 
eigenvalue $nH/2$ is a direct consequence of 
 (\ref{dira-extr2}). 

\qed


\begin{remark}
{\rm If $R$ denotes the scalar curvature 
of the induced metric on the embedded hypersurface $\Sigma$, 
we have the  Friedrich Inequality \cite{Fr}
\beq\label{lowe-boun-frie}
\lambda_1^2\ge \frac{n}{4(n-1)}\inf_\Sigma R.
\eeq
A consequence from the Gauss formula for the embedding 
$\Sigma\subset\Omega$ is that
$$
R=\overline{R}-2{\overline{\hbox{\rm Ric}}\,}(N,N)+n^2H^2-|\sigma|^2,
$$
where $\overline{\rm Ric}$ is the Ricci tensor of $\Omega$ and $\sigma$ is 
the second fundamental form of the embedding. From this 
equation it is clear that, in general, we cannot hope to 
obtain a relation between $R$ and $H$ allowing us to 
compare the Friedrich inequality (\ref{lowe-boun-frie}) and 
(\ref{lowe-boun-H}). 
However, 
if the Einstein tensor 
$\overline{\rm Ric}-(\overline{R}/2)\langle\;,\;\rangle$ of the ambient manifold 
is positive semidefinite, then
$$
R\le n^2H^2-|\sigma|^2\le n(n-1)H^2,
$$
where the last inequality is true because of the Schwarz 
inequality.  Moreover, if the last inequality is in fact an equality, then
the embedding is totally umbilical. Hence, we can state that}
\begin{quote}
When the Einstein tensor of the ambient manifold $\Omega$ 
is positive semidefinite, then the extrinsic lower bound (\ref{lowe-boun-H})
for the first eigenvalue of the Dirac operator of an embedded 
hypersurface  is sharper than the corresponding 
intrinsic Friedrich inequality (\ref{lowe-boun-frie}).
\end{quote}
{\rm This is the situation, for example, when the ambient 
space $M$ is Euclidean. In this case, if the scalar curvature 
$R$ of the hypersurface $\Sigma$ is positive, then $H$ 
should be positive too (see Lemma 1 in \cite{MR}), but it is 
possible for an embedded hypersurface  to have positive 
mean curvature (always with respect to the inner normal) and negative 
scalar curvature (for example, consider the 
revolution tori in ${\R}^3$). So there are situations in which only 
inequality (\ref{lowe-boun-H}) will be significant.}
\end{remark}

It is well known that lower bounds for the first 
eigenvalue of the Dirac operator on a compact Riemannian 
spin manifold (see \cite{Li}) 
have an important topological consequence: such a manifold with 
positive scalar curvature must have zero {\^A}-genus. This occurs 
because this topological invariant can be expressed in terms 
of the index of the Dirac operator and inequality (\ref{lowe-boun-frie}) 
imply that there are no non-trivial harmonic spinors. 
The same argument could be used to 
provide a topological 
obstruction for a compact hypersurface of 
a Riemannian spin manifold, with positive 
mean curvature  and nonnegative scalar 
curvature, to bound a domain.

\section{Constant mean curvature and minimal embedded hypersurfaces}

We shall assume now that there exists a non-trivial parallel spinor 
field $\psi_0$ on the Riemannian spin ambient manifold $M$. As we 
have pointed out before, this implies that $M$ is Ricci flat and 
reduces the possibilities for the (restricted) holonomy group of 
the manifold. In fact, when the manifold $M$ is complete, simply 
connected and irreducible, the existence of such a spinor field 
is equivalent (see \cite{W}) to one of the following facts: $M$ 
is either flat, a Calabi-Yau K{\"a}hler manifold (that is, its holonomy 
group is  $SU(m)$), a hyper-K{\"a}hler manifold (that is, its 
holomomy group is  $Sp(m)$), $\dim M=8$ and $M$ has holonomy group 
Spin$(7)$ or $\dim M=7$ and the holonomy group of $M$ is G$_2$. It 
is clear from (\ref{dira-extr2}) that one has 
$$
D\psi_0=\frac{nH}{2}\psi_0
$$
on the hypersurface $\Sigma$. Suppose also that the mean curvature $H$ 
of $\Sigma$ is constant and that this constant is nonnegative when 
computed with respect to the inner normal vector. Then we have  
$\lambda_1\le nH/2$ (this has been shown by B{\"a}r in \cite{Ba} when 
the ambient space is Euclidean). But, in the present situation, we can apply 
Theorem \ref{infH} and so we obtain that, in fact, $\lambda_1=nH/2$. 
Moreover each eigenspinor on $\Sigma$ associated with this first eigenvalue 
is the restriction to $\Sigma$ of a parallel spinor on the enclosed domain 
$\Omega$. From this, we can deduce the two results. First, we give
 a spinorial proof, in the spirit of the proof discovered 
by Reilly in \cite{Re}, for the celebrated Alexandrov Theorem 
(see \cite{A} and \cite{MR,O} for comments and references).

\begin{theorem}[The Alexandrov Theorem]\label{theo-alex}
The only compact embedded hypersurfaces in the Euclidean space with 
constant mean curvature are the round spheres.
\end{theorem}
\pf Let $\Sigma$ be such a hypersurface in ${\R}^{n+1}$ and 
$\Omega$ the enclosed domain. By studying the maxima on $\Sigma$ of 
the distance function  to a point in $\Omega$, one can see that the 
mean curvature $H$ of $\Sigma$ with respect to the inner normal must 
be positive (see \cite{MR} for more details). Then we are in the 
situation quoted before the statement of the theorem. That is, the 
first eigenvalue of the Dirac operator $D$ of the induced spin 
structure on $\Sigma$ is $nH/2$ and each corresponding eigenspinor 
is the restriction to $\Sigma$ of a parallel spinor on $\Omega$. If
$\xi(p)$ {\it denotes} the position vector at $p\in{\R}^{n+1}$, then
the spinor field $\psi$ defined, for all $p\in\Sigma$ by
$$
\psi(p)=\gamma\Bigl(H\xi(p)+N(p)\Bigr)\psi_0(p)$$
is also an eigenspinor on the hypersurface 
$\Sigma$ associated  with the first eigenvalue. For this, 
it is sufficient to consider (\ref{comp3}), identity 
(\ref{dira-extr1}) and the fact that
$$
\sum_{i=1}^n\gamma(e_i)\overline{\nabla}_{e_i}\psi=
\sum_{i=1}^{n}\gamma(e_i)\gamma(He_i-Ae_i)\psi_0=-n\hbox{\, trace\,} 
(HI-A)=0
$$
when $H$ is a constant. Hence this spinor field $\psi$ is the 
restriction to $\Sigma$ of a parallel spinor on $\Omega$. Then, 
for any vector field $X\in \Gamma(T\Sigma)$, we have
$$
0=\overline{\nabla}_X\psi=\gamma(H X-A X)\psi_0,$$
where $A$ is the shape operator corresponding to the normal 
field $N$. As $\psi_0$ has constant non-trivial length, we 
deduce that for any $X\in\Gamma(T\Sigma)$, $H X - AX=0$. 
Hence, the hypersurface is umbilical. As it is also compact, 
it is the round sphere.

\qed

An analogous reasoning works only for minimal hypersurfaces, 
if we replace the Euclidean ambient space by a Riemannian spin 
manifold with a non-trivial parallel spinor.

\begin{theorem}\label{theo-tota-geod}
Minimal compact hypersurfaces bounding a compact domain, in a 
Riemannian spin manifold admitting a non-trivial parallel 
spinor, are necessarily totally geodesic. This is the case for
minimal compact hypersurfaces embedded in simply connected 
Calabi-Yau manifolds, hyper-K{\"a}hler manifolds, or manifolds 
with holonomy groups {\rm Spin}$(7)$ and {\rm G}$_2$.
\end{theorem}
\pf Let $\psi_0$ be a non-trivial parallel spinor 
on the ambient manifold. As in the proof above, under these
hypotheses, we know that the first nonnegative eigenvalue  
of the Dirac operator $D$ on the hypersurface $\Sigma$ is 
$\lambda_1=0$, because $H=0$, and that the restriction of 
$\psi_0$ on $\Sigma$ is a corresponding eigenspinor. But 
 (\ref{anti-comm}) implies that 
also $D\Big(\gamma(N)\psi_0\Big)=0$, where $N$ is a unit normal 
field. Thus the spinor field $\gamma(N)\psi_0$ is the 
restriction to $\Sigma$ of a parallel spinor on the 
enclosed domain $\Omega$. Then, if $X\in \Gamma(T\Sigma)$, 
we have
$$
0=\overline{\nabla}_X\Big(\gamma(N)\psi_0\Big)=-\gamma(A X)\psi_0.
$$
Hence the shape operator $A$ is identically zero on 
$\Sigma$, that is, $\Sigma$ is totally geodesic.

\qed

Finally, we prove a rigidity result for 
compact hypersurfaces if the ambient Riemannian spin manifold 
admits a non-trivial parallel spinor. For this, we  
employ the same integral inequalities as above in a way 
inspired from the paper \cite{Ro} by A. Ros.

\begin{theorem}
Let $M$ be a Riemannian spin manifold with a 
non-trivial parallel spinor and $\Sigma\subset M$, a compact 
hypersurface with nonnegative mean curvature, bounding a compact 
domain in $M$. Suppose that $\iota :\Sigma\rightarrow M$ is 
another isometric immersion which induces on $\Sigma$ 
the same spin structure and whose mean curvature $H_\iota $ 
is constant and satis\-fies $|H_\iota |\le H$. Then the 
second fundamental form of the immersion $\iota $ 
coincides with that of the embedding $\Sigma\subset M$.
In particular, $|H_\iota |=H$.
\end{theorem}
\pf Under these assumptions, we can identify the 
spinor bundles induced on $\Sigma$ from the corresponding 
one on $M$, and the connections and Dirac operators 
constructed by means of the embedding $\Sigma\subset M$ 
and the immersion $\iota $. {}From a non-trivial parallel 
spinor on $M$, by using equation (\ref{dira-extr2}) for 
the immersion $\iota $, we get a spinor field $\varphi$ on 
$\Sigma$ such that
$$
D\varphi=\frac{n}{2}H_\iota \varphi\qquad\hbox{and }|\varphi|^2
=1\hbox{ on }\Sigma,
$$
where the mean curvature $H_\iota \ge 0$ is computed 
with respect to a suitable unit normal field $N_\iota $ 
for $\iota $ (note that $\Sigma$ is orientable because 
it bounds a domain). 
In particular, since $H_\iota $ is a nonnegative constant, 
we have $\pi_+\varphi=\varphi$. On the other hand, taking into 
account that $M$ is Ricci flat and that the mean curvature 
$H$ of the embedding $\Sigma\subset M$ is nonnegative, 
we can apply Proposition \ref{prop-APS-boun} to obtain a smooth solution of 
the following problem
\beq\left\{
\begin{array}{lll}
\overline{D}\psi &=0\qquad &\hbox{ {\rm on} }\Omega\\
\pi_+\psi&=\pi_+\varphi=\varphi\qquad &\hbox{ {\rm on} }\partial\Omega=\Sigma.
\end{array}
\right.
\eeq
Then, from inequalities (\ref{boun-weit-twis-ineq}) and 
(\ref{proj}), we deduce that
$$
0\le \int_\Sigma\Bigl(\langle D\pi_+\psi,\pi_+\psi\rangle
-\frac{nH}{2}|\psi|^2\Bigr)\,d\Sigma
$$
and the equality holds if and only if $\psi$ is parallel 
on $\Omega$ and $\psi=\pi_+\psi$. But we know that 
$$
D\pi_+\psi=D\varphi=\frac{n}{2}H_\iota \varphi
=\frac{n}{2}H_\iota \pi_+\psi.
$$
Hence
$$
0\le \int_\Sigma\Bigl(H_\iota  |\pi_+\psi|^2
-H|\psi|^2\Bigr)\,d\Sigma\le 0
$$
because $0\le H_\iota  \le H$ and $|\pi_+\psi|\le |\psi|$. 
As the equality is achieved and we know that 
$|\pi_+\psi|=|\varphi|=1$, we have  $\psi=\pi_+\psi=\varphi$,
and  $H_\iota =H$, hence the embedding $\Sigma\subset M$ 
has constant mean curvature $H$. Now, from (\ref{spin-gaus}) 
and the fact that $\psi$ is parallel on $\Omega$, we deduce 
$$
\gamma(A X)\gamma(N)\psi=\gamma(A_\iota  X)\gamma(N_\iota )\psi 
$$
for any $X\in \Gamma(T\Sigma)$, where $A_\iota $ is the shape 
operator corresponding to the orientation $N_\iota $ for the 
immersion $\iota $. Now, if $H=0$, then Theorem \ref{theo-tota-geod} above 
implies that $A$ is identically zero and hence $A_\iota\equiv 0$ too. 
Thus, assume  $H>0$. Multiplying both sides in 
the last equality by $\gamma(X)$ and contracting with $X$, 
since $H=H_\iota $, we have  $$
\gamma(N)\psi=\gamma(N_\iota )\psi.$$
As a consequence, if $X$ is tangent to the hypersurface 
$\Sigma$, we have 
$$
\gamma(A X)\psi=-\overline{\nabla}_X\Big(\gamma(N)\psi\Big)
=-\overline{\nabla}_X\Big(\gamma(N_\iota )\psi\Big)= 
\gamma(A_\iota  X)\psi.
$$
This implies that the two shape operators coincide, 
as claimed.

\qed

\end{document}